\documentclass[12pt,twoside,a4paper]{article}

\usepackage{amssymb}
\usepackage{amsfonts}
\usepackage{amsmath}
\newtheorem{example}{Example}[section]

\def\<{\langle}
\def\>{\rangle}
\def\R{{\mathbb R}}
\def\C{{\mathbb C}}
\def\Z{{\mathbb Z}}

\def\Hf{{\rm Hf\,}}
\def\SG{{\mathfrak S}}
\def\H{{\cal H}}

\def\ie{{\it i.e. }}

\def\goth{\mathfrak}

\def\shuff#1#2{\mathbin{
\hbox{\vbox{ \hbox{\vrule \hskip#2 \vrule height#1 width 0pt
}%
\hrule}%
\vbox{ \hbox{\vrule \hskip#2 \vrule height#1 width 0pt
\vrule }%
\hrule}%
}}}

\def\shuf{{\mathchoice{\shuff{7pt}{3.5pt}}%
{\shuff{6pt}{3pt}}%
{\shuff{4pt}{2pt}}%
{\shuff{3pt}{1.5pt}}}}%
\def\shuffle{\,\shuf\,}

\def\ashuff#1#2#3{
\kern 1pt \vrule height#1 \overline{\vrule height#3 width 0pt
\hskip#2} \rule{.3pt}{#1}\overline{\vrule height#3 width 0pt
\hskip#2} \rule{.3pt}{#1} \kern 1pt }

\def\ashuf{{\mathchoice{\ashuff{7pt}{3.5pt}{6pt}}%
{\ashuff{6pt}{3pt}{5pt}}%
{\ashuff{4pt}{2pt}{3.5pt}}%
{\ashuff{3pt}{1.5pt}{2.5pt}}}}%
\def\ashuffle{\,\ashuf\,}

\def\Pfa{{\rm Pf\,}}
\def\Haff{{\rm Hf\,}}

\def\ov#1{\overline{#1}}

\title{Pfaffian and Hafnian Identities in Shuffle Algebras}

\author{Jean-Gabriel {\sc Luque} and Jean-Yves {\sc Thibon}}

\date{}

\begin{document}

\maketitle

\begin{abstract}
Chen's lemma on iterated integrals implies that certain identities
involving multiple integrals, such as the de Bruijn and Wick
formulas, amount to combinatorial identities for Pfaffians and
hafnians in shuffle algebras. We provide direct algebraic proofs of
such shuffle identities, and obtain various generalizations. We
also discuss some Pfaffian identities due to Sundquist and
Ishikawa-Wakayama, and a Cauchy formula for anticommutative
symmetric functions. Finally, we extend some of the previous
considerations to hyperpfaffians and hyperhafnians.
\end{abstract}

\section{Introduction}

The free associative algebra $\Z\<A\>$ over an alphabet $A$ is
naturally endowed with a commutative operation, the shuffle
product $\shuffle$ which can be recursively defined by $$
au\shuffle bv = a(u\shuffle bv)+ b(au\shuffle v) $$ where $a,b\in
A$ and $u,v\in A^*$.

This operation has at least two interesting interpretations.
First, it is dual to the comultiplication  which admits the Lie
polynomials as primitive elements (Friedrich's theorem, see, {\it
e.g.},  \cite{Re}). Second, it describes the multiplication of
iterated integrals (Chen's lemma \cite{Chen2}). It is this last
property which will be of interest to us. To state it in a
convenient way, let us fix some vector space $\H$ of suitably
integrable functions over an interval $(a,b)$. For
$f_1,\ldots,f_n\in\H$, we identify the tensor product
$F=f_1\otimes\cdots\otimes f_n$ with the function of $n$ variables
$F(x_1,\ldots,x_n)=f_1(x_1)f_2(x_2)\cdots f_n(x_n)$. Let $\<\ \>$
be the linear form defined on each $\H^{\otimes n}$ by
\begin{equation}
\<f_1\otimes\cdots\otimes f_n\>= \int_a^{b}dx_1 \int_a^{x_1}dx_2
\cdots \int_a^{x_{n-1}}dx_n f_1(x_1)f_2(x_2)\cdots f_n(x_n)\,.
\end{equation}
We now choose a family of functions $(\phi_i)_{i\in I}$, labelled
by some alphabet $I$, and for a word $w=i_1i_2\cdots i_n\in I^*$,
we set
\begin{equation}
\<w\> =
\<\phi_{i_1}\otimes\phi_{i_2}\otimes\cdots\otimes\phi_{i_n}\>\,.
\end{equation}
Then, Chen's lemma states that
\begin{equation}\label{Chen}
\<u\>\<v\> = \<u\shuffle v\>\,.
\end{equation}

Iterated integrals occur for instance in the power series
solutions of linear differential equations with variable
coefficients. An example of current interest is provided by the
multiple zeta values \cite{Zag}, and an interesting application of
the shuffle product appears in  \cite{BBBL}, where a conjecture of
Zagier \cite{Zag}  is proved by reduction  to a simple shuffle
identity.

In certain cases, iterated integrals can be evaluated in closed
form by means of de Bruijn's formulas \cite{dB}

\begin{equation}\label{fdB1}
\mathop{\int\cdots\int}_{a\leq x_1< \cdots< x_{2n}\leq b} \det
\left(\phi_i(x_j)\right)dx_1\cdots
dx_{2n}=\Pfa\left(P_{ij}\right)_{1\leq i,j\leq 2n}
\end{equation}
where $\Pfa(P)$ denotes the Pfaffian of the antisymmetric
matrix $P$,
$$
P_{ij}=\mathop{\int\int}_{a\leq x<y\leq
b}[\phi_i(x)\phi_j(y)-\phi_j(x)\phi_i(y)]dxdy
$$
and
\begin{equation}\label{fdB2}
\mathop{\int\cdots\int}_{a\leq x_1< \cdots< x_{n}\leq b} \det
\left(\phi_i(x_j)|\psi_i(x_j)\right)dx_1\cdots
dx_n=\Pfa\left(Q_{ij}\right)_{1\leq i,j\leq 2n}
\end{equation}
where $Q_{ij}=\int_a^b[\phi_i(x)\psi_j(x)-\phi_j(x)\psi_i(x)]dx$
 and $\left(\phi_i(x_j)|\psi_i(x_j)\right)$ denotes the matrix
whose $i$th row is
$[\phi_i(x_1),\psi_i(x_1),\phi_i(x_2),\psi_i(x_2),\ldots,\phi_i(x_n),\psi_i(x_n)]$.

There is also a version of formula (\ref{fdB1}) for determinants
of odd order,
\[\begin{array}{ll}
\displaystyle \mathop{\int\cdots\int}_{a\leq x_1< \cdots<
x_{2n+1}\leq b} \det \left(\phi_i(x_j)\right)&dx_1\cdots
dx_{2n+1}\\
&=\displaystyle\sum_{p=1}^{2n+1}(-1)^{p+1}\int_a^b\phi_p(x)dx\
\Pfa\left(P_{ij }\right)_{1\leq i,j\leq 2n\atop i,j\neq p}\,.
\end{array}
\]
The foundations of the theory of random matrices can be developed
from these identities \cite{TW}.

\section{The Wick formulas}

The above formulas are reminescent of identities which are usually
called {\it Wick formulas} in the physics literature. These can be
stated in various ways, and come in two flavors: bosonic and
fermionic.

The bosonic version amounts to a familar fact about moments of
Gaussian measures. Let $d\mu(x)=C\, e^{-\frac12 {}^txAx}dx$ be a
Gaussian probability mesure on $\R^N$. We set $\langle
f\rangle=\int_{\R^N}f(x)d\mu(x)$. Let $f_1,\dots,f_i \dots$ be
linear forms. Then, the bosonic Wick formula reads
\begin{equation}\label{fbw}
\langle f_1\cdots f_n\rangle=\left\{\begin{array}{ll} 0&\mbox{if
}n \mbox{ is odd}\\ \Haff(\langle f_if_j\rangle)_{1\leq i,j\leq
n}&\mbox{if }n\mbox{ is even}\,,
\end{array}\right. 
\end{equation}
the {\it hafnian} of a symmetric matrix of even order $n=2r$ being
defined by
\begin{equation}
\Haff(A)=\sum_{I,J} a_{i_1j_1}a_{i_2j_2}  \cdots a_{i_rj_r}
\end{equation}
where the sum runs over all decompositions of $\{1,2,\ldots,2r\}$
into disjoints subsets $I=\{i_k\},\ J=\{j_k\}$, such that
$i_1<\cdots<i_r$, $j_1<\cdots <j_r$ and for each $k$, $i_k<j_k$.
Replacing $A$ by an antisymmetric matrix and multiplying each
term of the sum by the signature of the permutation $i_1j_1\cdots
i_rj_r$, one would obtain the Pfaffian of $A$. Actually, hafnians
were introduced by Caianiello \cite{Ca} in order  to emphasize the
similarities between the bosonic and fermionic versions of the
Wick formula.

These similarities become particularly obvious when the fermionic
Wick formula is stated in terms of {\it Berezin integrals}
(integration over Grassmann variables \cite{Be}).

The Berezin integral is essentially a convenient notation for
computing in Grassmann algebras. Let $\eta=\{\eta_i\}_{1\leq i\leq
N}$ be a finite set of anticommuting symbols (\ie
$\eta_i\eta_j+\eta_j\eta_i=0$). We denote by $\bigwedge_K(\eta)$
the $K$-algebra generated by the $\eta_i$. The Grassmann integral
of an anticommutative polynomial $f$ in the $\eta_i$ is defined by
\begin{equation}
\int d\eta_n\cdots d\eta_1
f(\eta_1,\dots,\eta_n)=\frac\partial{\partial\eta_n}\cdots
\frac\partial{\partial\eta_1}f(\eta_1, \cdots,\eta_n)
\end{equation}
where $\frac\partial{\partial\eta_i}$ acts on the Grassmann
algebra as a left derivation (we push $\eta_i$ to the left, with a
sign, and erase it). Therefore, in the fermionic calculus, 
integration and  derivation are one and the same operation.

Let $d\mu(\eta)=C e^{\frac12^t\eta Q\eta}$ be a ``fermionic
Gaussian probability mesure'' ($Q$ being an antisymetric and non
degenerate matrix, and $\int d\mu(\eta)=1$) and let
$f_1,\dots,f_n$ be linear forms in the $\eta_i$. We define
\begin{equation}
\langle f\rangle=\int d\mu(\eta) f(\eta)
\end{equation}
The fermionic Wick formula reads \cite{Be}
\begin{equation}\label{fwif}
\langle f_1\cdots f_n\rangle=\left\{
\begin{array}{ll}
0&\mbox{if }n\mbox{ is odd}\\ \Pfa\left(\langle f_i
f_j\rangle\right)_{1\leq i,j\leq n}&\mbox{if }n\mbox{ is even}\,.
\end{array}
\right.
\end{equation}

As we shall see, all these formulas are closely related, and
reflect simple combinatorial identities in the shuffle algebra.

\section{Pfaffians and Hafnians in  shuffle algebras}\label{s2}

\subsection{Wick formulas in shuffle algebras}

The Wick formula (\ref{fwif}) is equivalent to the equality
\begin{equation}\label{fwpdB}
\sum_{n\geq 0}\sum_{i_1<\cdots<i_n}\langle f_{i_1} \cdots
f_{i_{n}}\rangle
\eta_{i_1}\cdots\eta_{i_{n}}=\prod_{i}^\rightarrow\left(1+\sum_{i<j}
\langle f_{i} f_{j}\rangle\eta_{i}\eta_{j}\right)
\end{equation}

More generally, to an antisymmetric matrix $Q$ of  order $N$ over
a commutative ring $K$, we can associate the following element of
the Grassmann algebra $\bigwedge_K(\eta)$ over  $K$
\begin{equation}\label{serieferm}
T=\prod_i^\rightarrow\left(1+\sum_{i<j}Q_{i j}\eta_i\eta_j\right)
\end{equation}
and this product being commutative we can express it as an
exponential
\begin{equation}
T=\exp\left\{\frac12{^t\eta} Q\eta\right\}=\sum_{n\ge 0}
\frac1{n!}\left\{\frac12{^t\eta} Q\eta\right\}^n
\end{equation}
where $\eta$ denotes here the column vector $(\eta_i)$. The
coefficient  $\langle \eta_I|T\rangle$ of $\eta_I=\eta_{i_1}\cdots
\eta_{i_{2n}}$ in $T$ being $\displaystyle\int d\eta_{i_{2n}}
\cdots d\eta_{i_1}T$, the fermionic Wick formula gives
\begin{equation}\label{wickshuffle}
 \langle \eta_I|T\rangle=\Pfa\left(Q_{i_k i_l}\right)_{1\leq k,l\leq 2n}
 \end{equation}

Let us suppose now that $K=R\langle A\rangle_{\shuffle}$ is the
shuffle algebra over some commutative ring $R$. We will denote by
$\Pfa_{\shuffle}(M)$ (resp. $\Hf_{\shuffle}(M)$) the Pfaffian
(resp. hafnian) of a matrix $M\in R\langle A\rangle_{\shuffle}$.\\
Let $(a_i)$ and $(b_i)$ be two sequences
of  letters of $A$, and set $Q_{ij}=a_ib_j-a_jb_i$. Remarking
that
\begin{multline}
\sum_{j_1,\cdots, j_{k+1}} \left(a_{i_1}b_{j_1} \cdots
a_{i_{k}}b_{j_k}\shuffle a_{i_{k+1}}b_{j_{k+1}}\right)
\eta_{i_1}\eta_{j_1}\cdots \eta_{i_{k+1}}\eta_{j_{k+1}}
\\
=\sum_{j_1,\cdots, j_{k+1}}\sum_l a_{i_1} b_{j_1}\cdots
a_{i_{l-1}}b_{j_{l-1}}a_{i_{k+1}} b_{j_{k+1}}a_{i_l}b_{j_l}\cdots
a_{i_k}b_{j_k}\eta_{i_1}\eta_{j_1}\cdots
\eta_{i_{k+1}}\eta_{j_{k+1}}
\end{multline}
since the left-hand side is equal to the right-hand side, plus a
sum of terms which are symmetric in one pair of indices
$(j_r,j_s)$ and are therefore zero, we find
\begin{multline}\label{f2var}
\displaystyle\sum_{i_1<\cdots<i_{2k}}
\displaystyle\sum_{\sigma\in{\goth S}_{2k}} \epsilon(\sigma)
a_{i_{\sigma(1)}} b_{j_{\sigma(2)}}\cdots
a_{i_{\sigma(2k-1)}}b_{i_{\sigma(2k)}} \eta_{i_1}\cdots\eta_{
i_{2k}}\\ ={\displaystyle\prod_{i>0}^\rightarrow}
\left(1+\sum_{i<j} \left(a_ib_j-a_jb_ i \right)\eta_i\eta_j
\right)
\end{multline}
where $\epsilon(\sigma)$ is the signature of the permutation
$\sigma$, and formula (\ref{wickshuffle}) leads to
\begin{equation}\label{pfab}
\sum_{\sigma\in{\goth S}_{2n}}\epsilon(\sigma)
a_{{\sigma(1)}}b_{{\sigma(2)}}\cdots
a_{{\sigma(2n-1)}}b_{{\sigma(2n)}}
=\Pfa_{\shuffle}\left(a_ib_j-a_jb_i\right)_{1\leq
i,j\leq 2n}
\end{equation}
{\footnotesize
\begin{example}{\rm With $N=4$ we obtain
\[
\begin{array}{l}
a_1b_2a_3b_4-a_1b_2a_4b_3-a_1b_3a_2b_4+a_1b_4a_4b_2+
a_1b_4a_2b_3-a_1b_4a_3b_2\\
-a_2b_1a_3b_4+a_2a_4a_2b_3+a_2b_3a_1b_4 -a_2b_3a_4b_1+
-a_2b_4a_1b_4+ a_2b_4a_3b_1\\ -a_3b_2a_1b_4+
a_3b_2a_4b_1+a_3b_1a_2b_4 -a_3b_1a_4b_2 -a_3b_4a_1b_2+
a_3b_4a_2b_1\\ -a_4b_2a_3b_1+a_4b_2a_1b_3+ a_4b_3a_2b_1-
a_4b_3a_1b_2-a_4b_1a_2b_2+a_4b_1a_2b_3\\ =(a_1b_2-a_2b_1)\shuffle
(a_3b_3-a_4b_4)-(a_1b_3-a_3b_1)\shuffle(a_2b_4-a_4b_2)\\
+(a_1b_4-a_4b_1)\shuffle(a_2b_3-a_3b_2)

\end{array}
\]
}
\end{example}
} Similarly, if we consider the alphabet $A=\{a_{ij}\}$, we get
\begin{multline}\label{fshuffmat}
\displaystyle\sum_{i_1<\cdots<i_{2n}}\displaystyle\sum_{\sigma\in{\goth
S}_{2n}} \epsilon(\sigma)a_{i_{\sigma(1)}i_{\sigma(2)}}\cdots
a_{i_{\sigma(2n-1)}i_{\sigma(2n)}}\eta_{i_1}\cdots\eta_{ i_{2n}}\\
=\prod_i^\rightarrow\left(1+\sum_{i<j}\left(a_{ij}-
a_{ji}\right)\eta_i\eta_j\right)
\end{multline}
and the Wick formula gives
\begin{equation}\label{sdb2}
\sum_{\sigma\in{\goth S}_{2n}}\epsilon(\sigma)
a_{{\sigma(1)}{\sigma(2)}}\cdots a_{{\sigma(2n-1)}{\sigma(2n)}}=
\Pfa_{\shuffle}\left(a_{kl}-a_{lk}\right)_{1\leq k,l\leq 2n}
\end{equation}
Remark that although (\ref{pfab}) looks like the specialization
$a_{ij}=a_ib_j$ of (\ref{sdb2}), this is not the case, since this
specialization is not a shuffle homomorphism.

{\footnotesize
\begin{example}
{\rm with $1\leq i,j\leq 4$, the formula reads}
\begin{multline*}
(a_{12}-a_{21})\shuffle(a_{34}-a_{43})-(a_{13}-a_{31})\shuffle
(a_{24}-a_{42})+(a_{14}-a_{41})\shuffle (a_{23}-a_{32})\\
=a_{12}a_{34}-a_{12}a_{43}+a_{34}a_{12}-a_{43}a_{12}-a_{21}a_{34}
+a_{21}a_{43}-a_{34}a_{21}+a_{43}a_{21} \\
-a_{13}a_{24}+a_{13}a_{42}-a_{24}a_{13}+a_{42}a_{13}
+a_{31}a_{24}-a_{31}a_{42}+a_{24}a_{31}-a_{42}a_{31}\\
+a_{14}a_{32}-a_{14}a_{23}+a_{32}a_{14}-a_{23}a_{14}-
a_{41}a_{32}+a_{41}a_{23}-a_{32}a_{41}+a_{23}a_{41}.
\end{multline*}
\end{example}
}

Similar identities can be obtained by expanding a bosonic version
of the series (\ref{serieferm}). Let $\{\xi_i\}$ be a set of
commuting symbols verifying $\xi_i^2=0$, and let $Q$ be a
symmetric matrix with zero diagonal. The series
\begin{equation}
T=\prod_i\left(1+\sum_{i<j}Q_{ij}\xi_i\xi_j\right)
\end{equation}
has as coefficients
\begin{equation}
\langle \xi_I|T\rangle=\Hf\left(Q_{i_ki_l}\right)
\end{equation}
where $I=(i_1,\cdots,i_{2n})$. Specializing this to the shuffle
algebra over variables $a_{ij}$ and taking $Q_{kl}=a_{kl}+a_{lk}$
if $k\neq l$ and $Q_{kk}=0$, we obtain
\begin{equation}\label{fhaff2}
\sum_{\sigma\in{\goth S}_{2n}} a_{{\sigma(1)}{\sigma(2)}}\cdots
a_{{\sigma(2n-1)}{\sigma(2n)}}=\Haff_{\shuffle}
\left(Q_{kl}\right)_{1\leq k,l\leq 2n}\,.
\end{equation}

{\footnotesize
\begin{example}
{\rm Let $A=\{a_{ij}\}_{1\leq i,j\leq 4}$. We have}
\begin{multline*}
(a_{12}+a_{21})\shuffle(a_{34}+a_{43})+(a_{13}+a_{31})\shuffle
(a_{24}+a_{42})+(a_{14}+a_{41})\shuffle (a_{23}+a_{32})\\
=a_{12}a_{34}+a_{12}a_{43}+a_{34}a_{12}+a_{43}a_{12}+
a_{21}a_{34}+a_{21}a_{43}+a_{34}a_{21}+a_{43}a_{21}\\
+a_{13}a_{24}+a_{13}a_{42}+a_{24}a_{13}+a_{42}a_{13}+
a_{31}a_{24}+a_{31}a_{42}+a_{24}a_{31}+a_{42}a_{31}\\
+a_{14}a_{32}+a_{14}a_{23}+a_{32}a_{14}+a_{23}a_{14}+
a_{41}a_{32}+a_{41}a_{23}+a_{32}a_{41}+a_{23}a_{41}.
\end{multline*}
\end{example}
}

\subsection{Other Hafnian and Pfaffian identities}
The simplest one in the series is
\begin{equation}\label{fhaff1}
\sum_{\sigma\in{\goth S}_{2n}}a_{{\sigma(1)}}\cdots
a_{{\sigma(2n)}}=\frac1{(2n)!!}\Haff_{\shuffle}\left(Q_{kl}\right)_{1\leq
k,l\leq 2n}
\end{equation}
where $Q_{kl}=a_{k}a_{l}+a_{l}a_{k}$ if $k\neq l$ and $Q_{kk}=0$.
It amounts to the associativity and commutativity of the shuffle
product. {\footnotesize
\begin{example}
{\rm With $A=\{a_1,a_2,a_3,a_4\}$, we have
\begin{multline*}
\frac13\left[(a_1a_2+a_2a_1)\shuffle
(a_3a_4+a_4a_3)+(a_1a_3+a_3a_1)\shuffle(a_2a_4+a_4a_2)\right.\\
\left.+(a_1a_4+a_4a_1)\shuffle(a_2a_3+a_3a_2)\right] =
a_1a_2a_3a_4+a_1a_2a_4a_3+a_1a_3a_2a_4
+a_1a_3a_4a_2+\\a_1a_4a_2a_3+a_1a_4a_3a_2+
a_2a_1a_3a_4+a_2a_1a_4a_3+a_2a_3a_1a_4
+a_2a_3a_4a_1+a_2a_4a_1a_3+\\a_2a_4a_3a_1+
a_3a_2a_1a_4+a_3a_2a_4a_1+a_3a_1a_2a_4
+a_3a_1a_4a_2+a_3a_4a_1a_2+a_3a_4a_2a_1+\\
a_4a_2a_3a_1+a_4a_2a_1a_3+a_4a_3a_2a_1+
a_4a_3a_1a_2+a_4a_1a_3a_2+a_4 a_1a_2 a_3.
\end{multline*}
Indeed, the left hand side is just ${1\over 3}(3 \, a_1\shuffle
a_2 \shuffle a_3\shuffle a_4$). }
\end{example}
}

We can find other identities by considering the following
generating series in $\bigwedge_{K[x]}(\eta)$, where
$x=\{x_1,x_2,\cdots\}$ is a family of commuting variables:
\begin{equation}\label{fimpaire}
T(x,\eta)=\sum_{I\uparrow}t_I\eta_I=\prod_i^\rightarrow \left(
1+x_i\eta_i\right)
\end{equation}
where $I\uparrow$ means that $I$ is an increasing sequence. Since
$\eta_i^2=0$, one can write
\begin{equation}
T(x,\eta)=e^H=\prod^\rightarrow_ie^{x_i\eta_i}
\end{equation}
where
\begin{equation}
H=\sum_i x_i\eta_i+\frac12
\sum_{i<j}[x_i\eta_i,x_j\eta_j]=\sum_ix_i\eta_i+\sum_{i<j}x_ix_j\eta_i\eta_j
\end{equation}
is the Hausdorff series, here terminating in degree 2. On another
hand, one has
\begin{equation}\label{ftlgeneral}
t_I=\int d\eta_{\ov I}T(x,\eta)
\end{equation}
where $\ov I=(i_{n},i_{n-1},\cdots,i_1)$ if $I=(i_1,\cdots,i_n)$.
If $l( I)$ is even ($l(I)=2n$), we have as well
\begin{equation}\label{ftlpair}
t_I=\int d\eta_{\ov I}\
\exp\left\{\sum_{i<j}x_ix_j\eta_i\eta_j\right\}\,,
\end{equation}
and since
\begin{equation}
\exp\left\{\sum_{i<j}x_{i}x_{j}\eta_{i}\eta_{j}\right\}=\sum_n\sum_{i_1<j_1,
\dots,i_n<j_n \atop i_1<\cdots<i_n}x_{i_1}x_{j_1}\cdots
x_{i_n}x_{j_n}\eta_{i_1}\eta_{j_1}\cdots \eta_{i_n}\eta_{j_n}\,,
\end{equation}
we have finally
\begin{equation}\label{fimpairlpair}
t_I=\Pfa\left(Q_{kl}\right)_{1\leq k,l\leq 2n}
\end{equation}
where
\begin{equation}
Q_{kl}=\left\{\begin{array}{ll} x_{i_k}x_{i_l}&\mbox{ if } k<l\\
-x_{i_k}x_{i_l}&\mbox{ if } l<k\\ 0&\mbox{ if } l=k.
\end{array}\right.
\end{equation}
If $l( I)=2n+1$ is odd, (\ref{ftlgeneral}) gives
\begin{equation}\label{ftlimpair}
t_I=\int d\eta_{\ov I} \exp\left\{\sum_ix_i\eta_i\right\}
\exp\left\{\sum_{i<j}x_{i}x_{j}\eta_{i}\eta_{j}\right\}.
\end{equation}
But
\begin{equation}
e^{\sum_ix_i\eta_i}=1+\sum_ix_i\eta_i
\end{equation}
and then
\begin{equation}\label{fimpaitlimpair}
t_I=\sum_k(-1)^{k+1}x_{i_k}\Pfa(Q_{i j})_{1\leq i,j\leq 2n+1\atop
i,j\neq k}.
\end{equation}
Let us now specialize this formula to the shuffle algebra. In
$\bigwedge_{K\langle A\rangle_{\shuffle}}(\eta)$, the generating
series
\begin{equation}\label{ffimpairesh}
\sum_{I\uparrow}\sum_{\sigma\in {\goth
S}_n}a_{i_{\sigma(1)}}\cdots
a_{i_{\sigma(n)}}\eta_I=\prod_{i>0}^{\rightarrow}\left(1+a_i\eta_i\right)
\end{equation}
gives
\begin{equation}\label{wfimpaire}
\sum_{\sigma\in {\goth S}_n}a_{\sigma(1)}\cdots a_{\sigma(n)}=
\begin{cases}
\Pfa_{\shuffle}\left(Q_{kl}\right)_{1\leq k,l\leq n}&\text{if $n$
is even},\\ \sum_{p}(-1)^{p+1}a_{p}\shuffle
\Pfa_{\shuffle}\left(Q_{kl}\right)_{1\leq k,l\leq n\atop k,l\neq
p} &\text{if $n$ is odd},
\end{cases}
\end{equation}
where
\begin{equation}
Q_{k l}=\left\{
\begin{array}{ll}
a_{k}\shuffle a_{l}&\mbox{ if } k<l\\ -a_{k}\shuffle a_{l}&\mbox{
if } l<k\\ 0&\mbox{ if } k=l
\end{array}
\right.
\end{equation}
{\footnotesize
\begin{example}
{\rm For $A=\{a_1,a_2,a_3,a_4\}$, we have
\begin{multline*}
(a_1a_2+a_2a_1)\shuffle
(a_3a_4+a_4a_3)-(a_1a_3+a_3a_1)\shuffle(a_2a_4+a_4a_2)+\\
(a_1a_4+a_4a_1)\shuffle(a_2a_3+a_3a_2) =
a_1a_2a_3a_4+a_1a_2a_4a_3+ a_1a_3a_2a_4+a_1a_3a_4a_2+a_1a_4a_2a_3
+\\a_1a_4a_3a_2+a_2a_1a_3a_4+a_2a_1a_4a_3
+a_2a_3a_1a_4+a_2a_3a_4a_1+a_2a_4a_1a_3+
a_2a_4a_3a_1+a_3a_2a_1a_4+\\a_3a_2a_4a_1
+a_3a_1a_2a_4+a_3a_1a_4a_2+a_3a_4a_1a_2+
a_3a_4a_2a_1+a_4a_2a_3a_1+a_4a_2a_1a_3+
a_4a_3a_2a_1+\\a_4a_3a_1a_2+a_4a_1a_3a_2+a_4a_1a_2 a_3
\end{multline*}
Clearly, the left hand side reduces to $a_1\shuffle a_2\shuffle
a_3\shuffle a_4$, and the formula does not bring any new
information.}
\end{example}
} For completeness, let us remark that there is a bosonic version of
(\ref{fimpaire}) in $
K\langle A\rangle_{\ashuffle}[\xi]$, where $\ashuffle$ denotes the
antishuffle product on $K\langle A\rangle$ (i.e.
$\ashuffle=\shuffle_q$ for $q=-1$, where the $q$-shuffle is
defined by the recursive formula $au\shuffle_q bv =a(u\shuffle_q
bv)+q^{|au|}b(au\shuffle_q v)$)
\begin{equation}
{\prod_i^\rightarrow} \left(
1+x_i\xi_i\right)=\sum_{I\uparrow}\sum_{\sigma\in {\goth
S}_n}\epsilon(\sigma)a_{i_{\sigma(1)}}\cdots
a_{i_{\sigma(n)}}\xi_{I}
\end{equation}
leads to the following equality  in the antishuffle algebra
$K_{\ashuffle}\langle A\rangle$
\begin{equation}\label{fhaff3}
\sum_{\sigma\in{\goth S}_n}\epsilon(\sigma)a_{{\sigma(1)}}\cdots
a_{{\sigma(n)}}=
\begin{cases}
\Haff_{\ashuffle}\left(Q_{kl}\right)_{1\leq k,l\leq n}&\text{if
$n$ is even}\\
\sum_{p}a_{p}\ashuffle\Haff_{\ashuffle}\left(Q_{kl}\right)_{1\leq
k,l\leq n\atop k,l\neq p}&\text{if $n$ is odd}
\end{cases}
\end{equation}
where $Q_{k l}=Q_{l k}=a_{k}a_{l}-a_{l}a_{k}$ if $k\leq l$.

Expanding the hafnians, we see that this identity amounts to the
associativity and anticommutatitity of the antishuffle.

{\footnotesize
\begin{example}
{\rm Let A=$\{a_1,a_2,a_3,a_4\}$. We have
\begin{multline*}
(a_1a_2-a_2a_1)\ashuffle
(a_3a_4-a_4a_3)+(a_1a_3-a_3a_1)\ashuffle(a_2a_4-a_4a_2)+\\
(a_1a_4-a_4a_1)\ashuffle(a_2a_3-a_3a_2)
=
a_1a_2a_3a_4-a_1a_2a_4a_3-a_1a_3a_2a_4+
a_1a_3a_4a_2+a_1a_4a_2a_3+\\-a_1a_4a_3a_2-
a_2a_1a_3a_4+a_2a_1a_4a_3+a_2a_3a_1a_4+
-a_2a_3a_4a_1-a_2a_4a_1a_3+a_2a_4a_3a_1+\\
-a_3a_2a_1a_4+a_3a_2a_4a_1+a_3a_1a_2a_4
-a_3a_1a_4a_2-a_3a_4a_1a_2+a_3a_4a_2a_1-
a_4a_2a_3a_1+\\a_4a_2a_1a_3+a_4a_3a_2a_1+
-a_4a_3a_1a_2-a_4a_1a_3a_2+a_4a_1a_2 a_3.
\end{multline*}
Indeed, the left hand side reduces to $a_1\ashuffle a_2\ashuffle
a_3 \ashuffle a_4$. }
\end{example}
}

\subsection{Applications to iterated integrals}

Setting $a_i=b_i$ in (\ref{pfab}), we find
\begin{equation}\label{chen1}
\sum_{\sigma\in{\goth S}_{2n}}\epsilon(\sigma)
a_{{\sigma(1)}}a_{{\sigma(2)}}\cdots
a_{{\sigma(2n-1)}}a_{{\sigma(2n)}}
=\Pfa_{\shuffle}\left(a_{k}a_{l}-a_{l}a_{k}\right)_{1\leq k,l\leq
2n}
\end{equation}
Using Chen's lemma (\ref{Chen}), we recover  de Bruijn's formula (\ref{fdB1}).
Indeed, taking as alphabet $A=\{1,\ldots,2n\}$ the subscripts of
our functions $\phi_i$, and applying the linear form $\<\ \>$ to
both sides of (\ref{chen1}), we obtain
\begin{equation}
\sum_{\sigma\in{\goth S}_{2n}}\epsilon(\sigma)
\<\sigma\>=\<\Pfa_{\shuffle}\left(ij-ji\right)\>=\Pfa(\<ij-ji\>)
\end{equation}
which is exactly (\ref{fdB1}).

Formula (\ref{pfab}) can also be interpreted in terms of iterated
integrals, and the corresponding identity seems to be new. Taking
$4n$ functions $\phi_1,\ldots,\phi_{2n}$,
$\psi_1,\ldots,\psi_{2n}$ of $2n$ variables, we obtain
\begin{equation}
\mathop{\int\cdots\int}_{a\leq x_1< \cdots< x_{2n}\leq b} \det
\left(\phi_i(x_j)|\psi_i(x_j)\right)dx_1\cdots
dx_{2n}=\Pfa\left(Q_{ij}\right)_{1\leq i,j\leq 2n}
\end{equation}
where $(\phi_i(x_j)|\psi_i(x_j))$ denotes the $2n\times 2n$ matrix
whose $i$th row is
\[[\phi_i(x_1),\psi_i(x_2),\ldots,\phi_i(x_{2n-1}),\psi_i(x_{2n})],\]
and
\begin{equation}
Q_{ij}=\mathop{\int\int}_{a\leq x<y\leq
b}[\phi_i(x)\psi_j(y)-\phi_j(x)\psi_i(y)]dxdy\,.
\end{equation}
\footnotesize
\begin{example}\label{ex:su}\rm
Let us consider two families of complex numbers $(x_i)_{i\ge 1}$
and $(y_i)_{i\ge 1}$ (each $x_i$ being non zero). We construct for
each $i$ two functions
\[f_i(t)=\frac1{x^2_i}t^{\frac1{x_i}-1}\]
and
\[g_i(t)=\frac{y_i}{x_i}t^{\frac1{x_i}-1}.\]
The iterated integral \[ I_n=\mathop{\int\cdots\int}_{0\leq
t_1<t_2<\cdots<t_{2n}\leq 1}\det(f_i(t_j)|g_i(t_j))dt_1\cdots dt_n
\]
can be written as the Pfaffian
\[
I_n=\Pfa\left(\frac{y_j-y_i}{x_i+x_j}\right),
\]
The calculation of this Pfaffian, due to Sundquist \cite{Su}, is
discussed in section \ref{sec:examples}.
\end{example}
\normalsize

Application of Chen's lemma to (\ref{sdb2}) gives the second de
Bruijn formula (\ref{fdB2}), and we can deduce from (\ref{fhaff1})
a ``de Bruijn like'' equality for permanents
\begin{equation}\label{ifhaff1}
\begin{array}{ll}
\displaystyle\mathop{\int\cdots\int}_{a\leq t_1<\cdots<t_{2n}\leq
b} \mbox{\rm per}&\left(\phi_i(t_j)\right)dt_1\cdots dt_{2n}=\\ &
\left.
{\displaystyle\frac1{(2n)!!}}\right.\Haff\left(\displaystyle\mathop{\int\int
}_{a\leq x<y\leq
b}[\phi_i(x)\phi_j(y)+\phi_j(x)\phi_i(y)]dxdy\right).
\end{array}
\end{equation}
which is not very informative, as both sides are clearly equal to
the product of the integrals $\int_a^b\phi_i(x)dx$. 
The corresponding
identities for determinants and permanents of odd order can be
obtained in the same way.

\footnotesize
\begin{example}{\rm
We can write down a permanental version of Wigner's integral
\cite{Me}, also discussed in section \ref{sec:examples}
\[I_n=\mathop{\int\cdots\int}_{a\leq t_1<t_2<\cdots<t_{2n}\leq b}{\rm
per}\left(e^{x_it_j} \right)dt_1\cdots dt_{2n}\] which evaluates
to
\[
I_n
=
{1\over (2n)!!}\Haff\left(\frac{(b_i-a_i)(b_j-a_j)}{x_ix_j}\right)
=
\prod_i\frac{b_i-a _i}{x_i}
\] 
with $a_i=e^{ax_i}$ and $b_i=e^{bx_i}$. }
\end{example}

\begin{example}{\rm
The equality
\begin{equation}\label{sum1}
\sum_{\sigma\in\SG_m}\frac1{x_{\sigma(1)}}\frac1{x_{\sigma(1)}+
x_{\sigma(2)}}\cdots\frac1{x_{\sigma(1)}+\cdots+x_{\sigma(m)}}
=\frac1{x_1\cdots x_m}
\end{equation}
can be easily shown by induction on $m$, remarking that
\begin{equation}\label{sum2}
\sum_k \frac1{x_1\cdots x_{k-1}x_{k+1}\cdots
x_m}=\frac{x_1+\cdots+x_m}{x_1\cdots x_m}.
\end{equation}
But we can also see it as an integral of type (\ref{ifhaff1}).
Indeed, direct calculation of
\[
I(x_1,\cdots,x_{m})=\mathop{\int\cdots\int}_{0\leq t_1<\cdots<t_m\leq 1}
\mbox{\rm per}\left(t_i^{x_j-1}\right)dt_1\cdots dt_{m} 
\] 
 gives
\[
I(x_1,\cdots,x_{m})=\sum_{\sigma}\frac1{x_{\sigma(1)}}\frac1{x_{\sigma(1)
}+
x_{\sigma(2)}}\cdots\frac1{x_{\sigma(1)}+\cdots+x_{\sigma(m)}}\,.
\]
On another hand, 
\[
I(x_1,\cdots,x_{m})=\lambda_1\cdots\lambda_m
\] with
$\lambda_i=\int_0^1t^{x_i-1}dt=\frac1{x_i}$. 
}
\end{example}
\normalsize

Finally,  applying Chen's lemma to (\ref{fhaff2}) we obtain
\begin{multline}\label{ifhaff2}
\displaystyle\mathop{\int\cdots\int}_{a\leq t_1<\cdots<t_n\leq b}
{\rm per}\left(\phi_i(t_j)|\psi_i(t_j)\right)dt_1\cdots dt_n\\
=\Haff\left(\int_a^b[\phi_i(x)\psi_j(x)+\phi_j(x)\psi_i(x)]\right)
\end{multline}

\footnotesize
\begin{example}{\rm
The hafnian analog of example \ref{ex:su} is
\[
I(x,y)=\mathop{\int\cdots\int}_{0\leq t_1<\cdots <t_n\le
1}\mbox{\rm per}\left(y_it_j^{x_i-1}\right.\left|t^{x_i}_j\right)
dt_1\cdots dt_n
=
\Haff\left(\frac{y_i+y_j}{x_i+x_j}\right).
\]
If we compute $I(x_1,\cdots,x_{2n},y_1,\cdots,y_{2n})$ directly,
we find the symmetrization identity
\[
\sum_{\sigma\in{\goth S}_{2n}}\frac{y_{\sigma(1)}
y_{\sigma(3)}\cdots
y_{\sigma(2n-1)}}{(x_{\sigma(1)}+x_{\sigma(2)})
(x_{\sigma(1)}+x_{\sigma(2)}+x_{\sigma(3)}+x_{\sigma(4)})
\cdots(x_{\sigma(1)}+\cdots+ x_{\sigma(2n)})}
=\Haff\left(\frac{y_i+y_j}{x_i+x_j}\right).
\]}
\end{example}
\normalsize

\section{Examples}\label{sec:examples}

\def\Pf{{\rm Pf\,}}

\subsection{Pfaffian identities}

The classical identity of Schur \cite{Sch}
(see \cite{Kn} for an illuminating discussion
of this identity and of similar ones)
\begin{equation}\label{Schu}
\Pf\left({x_k-x_l\over x_k+x_l}\right)=
\prod_{k<l}\left({x_k-x_l\over x_k+x_l}\right)
\end{equation}
(where $k,l=1,\ldots,2m$) has been recently generalized by
Sundquist \cite{Su}, giving evaluations of $$
\Pf\left({y_k-y_l\over x_k+x_l}\right) \qquad {\rm and}\qquad
\Pf\left({y_k-y_l\over 1+x_kx_l}\right)\, $$ and also by Ishikawa
and Wakayama \cite{IW2}, giving the interpolation $$
\Pf\left({y_k-y_l\over a+ b(x_k+x_l)+cx_kx_l}\right)\,. $$
Actually, all three identities are equivalent. Indeed, one can
obtain $\Pf\left({y_k-y_l\over 1+x_kx_l}\right)$ from
$\Pf\left({y_k-y_l\over z_k+z_l}\right)$ by setting
$z_k=\frac12\left({1-ix_k\over 1+ix_k}\right)$ and then setting
$x_k=\alpha+\beta t_k$ one obtains $\Pf\left({y_k-y_l\over
(1+\alpha^2)+\alpha\beta(t_k+t_l)+\beta^2t_kt_l}\right)$.
Therefore, to relate these identities to the previous
considerations, it will be sufficient to consider one of them, say,
the one for $\Pf\left({y_k-y_l\over x_k+x_l}\right) $. To this
aim, we consider the  de Bruijn integral
\begin{equation}\label{wigner} I_{2m}= \mathop{\int\cdots\int}_{a\leq t_1<
\cdots< t_{2m}\leq b} \det \left[e^{x_kt_l}\right] dt_1\cdots
dt_{2m} = \Pf(Q)
\end{equation}
where $Q=(Q_{kl})$ is given in terms of $a_k=e^{ax_k}$ and
$b_l=e^{bx_l}$ by
\begin{equation}
Q_{kl}={1\over x_kx_l(x_k+x_l)} \left|
\begin{matrix}u_k & v_k \\ u_l & v_l \end{matrix}\right|
\end{equation}
and $u_k=b_k-a_k$,  $v_k=x_k(a_k+b_k)$. The integral
(\ref{wigner}), which is of interest in random matrix theory, can
be evaluated directly, without the help of de Bruijn's formula.
This is a rather tedious calculation, originally due to E. Wigner
(unpublished), which can be found in Mehta's book \cite{Me} (see
p. 455, A.23). This calculation results into an expression
different from the previous one (\cite{Me}, Theorem 10.9.1), which
in the present notation reads
\begin{equation}\label{mehta}
I_{2m}= \prod_{k=1}^{2m}{1\over x_k}\prod_{k<l}{1\over x_k+x_l}
\det(x_i^{j-1}(b_i+(-1)^ja_i)) \,.
\end{equation}
Comparing the right hand sides of (\ref{mehta}) and
(\ref{wigner}), and substituting $b_k-a_k=u_k$,
$b_k+a_k=x_k^{-1}v_k$, one finds
\begin{multline}\label{sund}
\begin{vmatrix}
u_1 & v_1 & x_1^2u_1 & x_1^2v_1 & x_1^4u_1 & \cdots &
x_1^{2m-2}v_1 \\ u_2 & v_2 & x_2^2u_2 & x_2^2v_2 & x_2^4u_2 &
\cdots & x_2^{2m-2}v_2 \\ \vdots&\vdots&\vdots & \vdots & \vdots&
\ddots &\vdots \\ u_{2m} & v_{2m} & x_{2m}^2u_{2m} &
x_{2m}^2v_{2m} & x_{2m}^4u_{2m} & \cdots & x_{2m}^{2m-2}v_{2m} \\
\end{vmatrix}
\\
=\prod_{i<j}(x_i+x_j)\Pf\left({u_iv_j-u_jv_i\over x_i+x_j}\right)
\end{multline}
which is one of Sundquist's formulas (when one sets $u_i=1$ and
$v_i=y_i$).

In Section \ref{s2}, we have interpreted de Bruijn's formula as a
combinatorial property of shuffles. In the following, we will
observe that (at least) parts of Wigner's argument can be
interpreted similarly.

To establish (\ref{mehta}), Wigner starts with the simplest
integral
\begin{equation}
J_n(x_1,\ldots,x_n) = \mathop{\int\cdots\int}_{a\leq t_1< \cdots<
t_{2m}\leq b} e^{t_1x_1+\cdots+t_nx_n}dt_1\cdots dt_n\,,
\end{equation}
 which yields $I_n$ by antisymmetrization:
\begin{equation}
I_n(x_1,\ldots,x_n)={\cal A}J_n(x_1,\ldots,x_n)
=\sum_{\sigma\in\SG_n}\varepsilon(\sigma)
J_n(x_{\sigma(1)},\ldots,x_{\sigma(n)})\,.
\end{equation}
If one computes the first integrals $J_1,J_2,\ldots$, an induction
pattern emerges, and one arrives at the expression
\begin{equation}\label{eqJ}
\begin{array}{l}
J_n(x_1,\ldots,x_n)=\\ \displaystyle\sum_{k=0}^n(-1)^k a_k\cdots
a_2a_1\cdot b_{k+1}\cdots b_n
R(x_k,x_{k-1},\ldots,x_1)R(x_{k+1},x_{k+2},\ldots,x_n)
\end{array}
\end{equation}
where we have set
\begin{equation}
R(z_1,\ldots,z_r)={1\over z_1(z_1+z_2)\cdots(z_1+z_2+\cdots
z_r)}\,.
\end{equation}

The rest of the calculation relies upon the two identities
\begin{equation}\label{mehta1}
\sum_{k=0}^n(-1)^k
R(x_k,x_{k-1},\ldots,x_1)R(x_{k+1},x_{k+2},\ldots,x_n) = 0
\end{equation}
and
\begin{equation}\label{mehta2}
{\cal A}R(x_1,\ldots,x_n)= \prod_{i=1}^{n}{1\over
x_i}\prod_{i<j}{x_j-x_i\over x_j+x_i} \,.
\end{equation}
Remark that if we set $\phi_z(t)=t^{z-1}$, we can write
\begin{equation}
R(z_1,\ldots,z_r)=\<z_1\cdots z_r\>
=\mathop{\int\cdots\int}_{0\leq t_1< \cdots< t_{r}\leq 1}
\phi_{z_1}(t_1)\cdots \phi_{z_r}(t_r) dt_1\cdots dt_r
\end{equation}
so that for $u=x_k\cdots x_1$, $v=x_{k+1}\cdots x_n$ (regarded as
words), $R(u)R(v)=\<u\>\<v\>=\<u\shuffle v\>$.

With this at hand, it is clear that (\ref{mehta1}) amounts to the
fact that the map $S(w)=(-1)^{|w|}\bar{w}$, where $\bar w$ is the
mirror image of the word $w$, is the antipode of the shuffle Hopf
algebra. This means that in the convolution algebra of
$K_{\shuffle} \<A\>$, i.e., ${\rm End^{\rm gr}}(K_{\shuffle}
\<A\>)$ endowed with the product $f*g={\rm sh}\circ (f\otimes
g)\circ C$, where ${\rm sh}(u\otimes v)=u\shuffle v$ and
$C(w)=\sum_{uv=w} u\otimes v$, one has $$ S*I=I*S=\iota\circ
\epsilon $$ where $\iota$ is the unit and $\epsilon$ the counit,
so that $S*I$ of a polynomial is equal to its constant term.
Hence, for any non empty word $w$, $$ \sum_{uv=w}(-1)^{|u|}\bar
u\shuffle v=0\,, $$ which yields (\ref{mehta1}). Equation
(\ref{mehta2}) can be obtained by writing $$ {\cal A}
R(x_1,\ldots,x_n)= \mathop{\int\cdots\int}_{0\leq t_1< \cdots<
t_{n}\leq 1} \det(\phi_{x_i}(t_j))dt_1\ldots dt_n = \Pf(A) $$
where $$A_{ij}=\int_0^1dt_j\int_0^{t_j}dt_i
(t_i^{x_i-1}t_j^{x_j-1}-t_j^{x_i-1}t_i^{x_j-1}) ={x_j-x_i\over
x_ix_j(x_i+x_j)} $$ and applying Schur's identity (\ref{Schu}).\\

There is also an interesting interpretation of (\ref{sund}),
involving the hyperoctahedral group. Let us introduce a family of
functions labelled by the alphabet $A=\{\bar n,\ldots,\bar 2,\bar
1,1,2,\ldots,n\}$
\begin{equation}
\phi_k(t)=b_kt^{x_k-1} \quad {\rm and }\quad \phi_{\bar
k}(t)=-a_kt^{x_k-1}\,,
\end{equation}
so that (\ref{eqJ}) can be written as
\begin{equation}
J_n(x_1,\ldots,x_n)= \left\< \sum_{k=0}^{n} \bar k \cdots \bar
2\bar 1\shuffle k+1\cdots n\right\>
\end{equation}
and
\begin{equation}\label{eqAJ}
{\cal A}J_n=\left\< \sum_{\sigma\in\SG_n}\epsilon(\sigma)\sigma
\sum_{k=0}^{n} \bar k \cdots \bar 2\bar 1\shuffle k+1\cdots
n\right\>
\end{equation}

If the word $w$ occurs in $\bar k \cdots \bar 2\bar 1\shuffle
k+1\cdots n$ with the letters $\bar k,\ldots,\bar 1$ at positions
$i_1,\ldots,i_k$, the sign of the underlying permutation $\alpha$,
obtained by erasing the bars, is $\epsilon(\alpha)=(-1)^{||I||}$,
where $||I||=\sum_l (i_l-1)$. We can furthermore interpret such a
word $w$ as an element of the hyperoctahedral group $B_n$,
generated by permutations and by the involutions $\tau_i :
i\mapsto \bar i$ and $\tau_i(j)=j$ for $j\not=i$. The double sum
in (\ref{eqAJ}) can now be interpreted as an element of the group
algebra $\Z B_n$, and since ${\cal A}\alpha=\epsilon(\alpha){\cal
A}$, we have for $w$ as above
\begin{equation}
{\cal A}w=(-1)^{||I||}{\cal A}\tau_{i_1}\cdots \tau_{i_k}\,.
\end{equation}
If we extend the action of $B_n$ to polynomials in
$x_i,y_i,a_i,b_i$ by setting $\bar x_i=x_i$, $\bar y_i=y_i$ and
$\bar b_i=-a_i$, (\ref{eqAJ}) can now be rewritten
\begin{equation}
{\cal A}J_n= {\cal A}\circ(1+\tau_1)\circ
(1-\tau_2)\circ\cdots\circ(1-\tau_{2m}) {b_1b_2\cdots b_{2m}\over
x_1(x_1+x_2)\cdots (x_1+\cdots+ x_{2m})}\,,
\end{equation}
that is,
\begin{equation}
{\cal A}J_n = {\cal A}(QR)
\end{equation}
where
\begin{equation}
Q = \prod_{i=1}^{2m}(b_i+(-1)^i a_i)
\end{equation}
%
Let $x^\rho=x_1^0x_2^1\cdots x_n^{n-1}$, so that $\Delta={\cal
A}(x^\rho)$. Then, we see that (\ref{sund}) is equivalent to the
symmetrization identity
\begin{equation}
{\cal A}(RQ){\cal A}(x^\rho)={\cal A}(R){\cal A}(Qx^\rho)\,.
\end{equation}

\subsection{Anticommutative symmetric functions and alterning quasi-symmetric
functions}

The fermionic Wick formula admits an  interpretation in the
language of non-commutative symmetric functions \cite{ncsf}. Let
us recall that the starting point of this theory is the
observation that the classical algebra of symmetric functions is a
polynomial algebra $Sym=\C[h_1,h_2\cdots]$ in an infinite number
of indeterminates (here, the $h_i$ are the complete homogeneous
symmetric functions). One  defines the algebra {\bf Sym} of
noncommutative symmetric functions as the free associative algebra
$\C\langle S_1,S_2,\cdots\rangle$ over an infinite sequence of
letters, together with a ``commutative image'' homomorphism
$S_n\rightarrow h_n$. In other words, {\bf Sym} is the tensor
algebra of the vector space $V=\bigoplus\C S_n$ and $Sym$ is its
symmetric algebra.

Here, we propose to have a look at its exterior algebra, which we
will call $ASym$ (anticommutative symmetric functions). It is the
image of {\bf Sym} under the specialization $S_n\rightarrow
S_n(A)=\eta_n$ (that is, the letter $A$ is used here as a label
for the image $F(A)$ of any  noncommutative symmetric function $F$
by the ring homomorphism $S_n\mapsto \eta_n$).

The dimension of the homogeneous component of degree $n$ of $ASym$
is equal to the number of partitions of $n$ into distinct parts.
As an algebra, $ASym$ is rather trivial, but the dual Hopf algebra
looks more interesting. To study it, we can follow the same
strategy as in \cite{ncsf}, and start with the Cauchy kernel. In
the notation of \cite{ncsf}
\begin{equation}
\begin{array}{rcl}
\sigma_1(XA)&=&\displaystyle\prod_{k\geq 1}^\rightarrow
\sigma_{x_k}(A)\,=\,\prod_{k\geq
1}^\rightarrow\sum_{i_k\geq0}x^{i_k}_k\eta_{i_k}\\
&=&\displaystyle\sum_{r}\sum_{i_1<\cdots<i_r}\left(\sum_{\sigma\in{\goth
S}_r}\epsilon(\sigma)M_{I\sigma}(X)\right)\eta_{i_1}\cdots\eta_{i_r}\\
&=&
\displaystyle\sum_{r}\sum_{i_1<\cdots<i_r}\left(\sum_{\sigma\in{\goth
S}_r}\epsilon(\sigma)M_{I\sigma}(X)\right)S^I(A)\end{array}
\end{equation}
where the $M_I(X)$ are the quasimonomial quasi-symmetric functions,
and $I \sigma$ means $(i_{\sigma(1)},\cdots,i_{\sigma(r)})$. The
coefficients
\begin{equation}
V_I(X)=\sum_{\sigma}\epsilon(\sigma)M_{I\sigma}(X)=\int
d\eta_{\bar I}\sigma_1(XA)
\end{equation}
 can be evaluated by means of the fermionic
Wick formula. Introducing the operators $$
\phi_k={\partial\over\partial\eta_k}=\int d\eta_k\   \quad
\mbox{\rm and}\ \phi^*_k(f)=\eta_kf\, $$ one has
\begin{equation}
V_I=\langle0|\phi_{i_r}\cdots\phi_{i_1}|\sigma(XA)\rangle
\end{equation}
where $|0\rangle=1$ and
$\displaystyle\sigma_1(XA)=\prod_{i\geq1}^\rightarrow
e^{\Phi(x_i)}|0\rangle$ if $\sigma_x(A)=e^{\Phi(x)}$ (here the
generating series $\Phi(x)$ of the second kind power sums reduces
to  $x\eta_1+x^2\eta_2+\cdots$). The Hausdorff series for the
product of exponentials is now
\begin{multline}
H=H(\Phi(x_1),\Phi(x_2),\cdots)\\ =\sum_{k\geq
1}\psi(x_k)+\frac12\sum_{i<j}[\Phi(x_i),\Phi(x_j)]\\ =\sum_{n\geq
1}M_n(X)\phi_n^*+\frac12\sum_{k\neq
l}\left(M_{kl}(X)-M_{lk}(X)\right)\phi_k^*\phi_l^*\\
=\frac12\sum_{k\neq l}Q_{kl}\phi^*_k\phi^*_l+\sum_{n\geq
1}M_n\phi^*_n
\end{multline}
where $Q_{kl}=M_{kl}-M_{lk}$. Hence,
\begin{equation}
V_I(X)=\int d\eta_{\bar I}e^H=
\langle0|e^{\sum_{n\geq0}M_n\phi^*_n}e^{{\frac12}^t\phi^*Q\phi^*}|0\rangle
\end{equation}
since the quadratic terms are central, so, expanding the
exponential $e^{\sum_nM_n\phi^*_n}$, we obtain the following
expression for the ``alternating quasisymmetric functions''.
\begin{multline}
V_I(X)=\sum_\sigma\epsilon(\sigma) M_{I\sigma}(X)\,
=\,\langle\phi_{i_r}\cdots\phi_{i_1}\rangle\\
=
\begin{cases}
\Pfa(Q_{i_ki_l})_{1\leq k,l\leq 2m}&\text{if $r=2m$ is even}\\
\sum_{k=1}^{2m+1}(-1)^{k-1}M_{i_k}\langle
\phi_{i_{2m+1}}\cdots\widehat{\phi_{i_k}}\cdots\phi_{i_1}\rangle&
\text{if $r =2m+1$ is odd}.
\end{cases}
\end{multline}

This formula is similar to Schur's expression of his symmetric
functions $Q_\lambda$ as Pfaffians of the $Q_{(i,j)}$. Remark that
these symmetric functions are also labelled by strict partitions,
and that Schur's Pfaffian can also be derived from the fermionic
Wick formula \cite{You}.

\section{Generalizations}
\subsection{Hyperpfaffians}
The Pfaffian of a skew symmetric  matrix of even order can be
expanded as
\begin{equation}
\Pfa\left(M_{ij}\right)_{1\leq i,j\leq 2n}=\sum_{\sigma\in{\goth
E}_{2n,2}}\epsilon(\sigma)M_{\sigma(1)\sigma(2)}\cdots
M_{\sigma(2n-1)\sigma(2n)}
\end{equation}
where the sum is over

\begin{multline}
{\goth E}_{2n,2}=\{\sigma\in{\goth
S}_{2n}|\sigma(2i+1)<\sigma(2i+2), \sigma(2p-1)<\sigma(2p+1),\\
0\leq i\leq n-1, 1\leq p\leq n-1\}
\end{multline}
We can extend this definition by setting
\begin{eqnarray*}
{\goth E}_{kn,k} &=&\{\sigma\in{\goth S}_{kn}|
\sigma(ki+1)<\cdots<\sigma(k(i+1)),\\
&&\sigma(k(p-1)+1)<\sigma(kp+1), 0\leq i\leq n-1, 1\leq p\leq
n-1\}.
\end{eqnarray*}
If we consider variables $M_{i_1\cdots i_k}$ verifying
$M_{i_{\sigma(1)}\cdots
i_{\sigma(k)}}=\epsilon(\sigma)M_{i_1\cdots i_k}$, the polynomial
\begin{equation}
{\Pfa}^{[k]}\left(M_{i_1\cdots i_k}\right)_{1\leq i_1,\dots,
i_k\leq nk}= \sum_{\sigma\in {\goth E}_{kn,k}}\epsilon(\sigma)
M_{{\sigma(1)}\dots {\sigma(k)}}\cdots M_{{\sigma((n-1)k+1)}\dots
{\sigma(nk)}}
\end{equation}
appears as a natural generalization of the Pfaffian. 
This polynomial has been called
the {\it hyperpfaffian} of the alternating tensor $M$ by Barvinok \cite{Ba}.

Similarly, higher order Hafnians, or {\it hyperhafnians}, can be
defined by
\begin{equation}
{\Haff}^{[k]}(M_{i_1\cdots i_k})_{1\leq i_1,\cdots,i_k\leq nk}=
\sum_{\sigma\in {\goth E}_{kn,k}} M_{{\sigma(1)}\dots
{\sigma(k)}}\cdots M_{{\sigma((n-1)k+1)}\dots {\sigma(nk)}}.
\end{equation}

Alternatively, these definitions can be presented as follows.
Given an alternating tensor $M$ of order $k$ as above, consider
the element
\begin{equation}
\Omega_M=\sum_{1\le i_1<\ldots <i_k\le nk} M_{i_1\cdots i_k}
\eta_{i_1}\cdots \eta_{i_k}
\end{equation}
of the Grassmann algebra of rank $nk$. Then, in this algebra
\begin{equation}
\Omega_M^n=n!{\Pfa}^{[k]}(M)\eta_1\eta_2\cdots \eta_{kn} \,.
\end{equation}
When $k=2$, this is one of the classical definitions of the
Pfaffian.

Similarly, when $M$ is a symmetric tensor of order $k$, again over
an $nk$ dimensional space, consider the element
\begin{equation}
G_M=\sum_{1\le i_1<\ldots <i_k\le nk} M_{i_1\cdots i_k}
\xi_{i_1}\cdots \xi_{i_k}
\end{equation}
wher $\xi_1,\ldots,\xi_{kn}$ are commuting symbols satisfying
$\xi_i^2=0$. Then
\begin{equation}
G_M^n=n! {\Haff}^{[k]}(M) \xi_1\xi_2\cdots \xi_{kn}\,.
\end{equation}

\subsection{Composition of hyperpfaffians}

The aim of this subsection is to compute the hyperfaffian
\[\Pf^{[2m]}\left(\Pf(a_{i,j})_{i,j\in\{k_1,\cdots,k_{2m}\}}\right)_{k_1\cdots,
k_{2m}\in \{1,\cdots, 2mn\}}\] where $A=(a_{i,j})_{1\leq i,j\leq
2mn}$ is a skew symmetric matrix and $m,n>0$. In order to
express this as a classical Pfaffian, we consider the $2-$form
$$
\omega=\sum_{1\leq i,j\leq2mn} a_{i j}\eta_i\eta_j \,. 
$$ 
The
$mn{\mbox{th}}$-power of $\omega$ is by definition of Pfaffians
$$
\omega^{nm}=(mn)!\Pf(a_{i,j})_{1\leq i,j\leq 2mn}\,. 
$$ 
But, on
another hand 
$$ 
\omega^{mn}=(m!\Omega)^n
$$ 
where
$$
\Omega=\sum_{1\leq
k_1<k_2<\cdots<k_{2m}}\Pf(a_{ij})_{i,j\in\{k_1,\cdots,k_{2m}\}}\eta_{k_1}
\cdots\eta_{k_{2m}}\,.
$$
Remarking that $\Omega^n$ is given by a hyperpfaffian 
$$
\Omega^n=n!\Pf^{[2m]}\left(\Pf(a_{ij})_{i,j\in\{k_1,\cdots,k_{2m}\}}\right)
_{1\leq
k_1\cdots, k_{2m}\leq 2mn}\eta_{1}\cdots\eta_{2nm}\,,$$ 
by
identification, we find
\begin{equation}
\Pf^{[2m]}\left(\Pf(a_{i,j})_{i,j\in\{k_1,\cdots,k_{2m}\}}\right)_{1\leq
k_1\cdots, k_{2m}\leq
2mn}=\frac{(mn!)}{(m!)^nn!}\Pf(a_{ij})_{1\leq i,j\leq 2mn}\,.
\end{equation}

\footnotesize
\begin{example}{\rm
If we set $m=2$ and $a_{i,j}=\frac{x_i-x_j}{x_i+x_j}$, we can
generalize  Schur's formula to the $4$-hyperpfaffian
\begin{equation}
\Pf^{[4]}\left(\prod_{1\leq s,t\leq
4}\frac{x_{k_s}-x_{k_t}}{x_{k_s}+x_{k_t}}\right)_{1\leq
k_1,k_2,k_3,k_4\leq 4n}=(2n)!!\prod_{1\leq i,j\leq
4n}\frac{x_i-x_j}{x_i+x_j}
\end{equation}}
\end{example}
\normalsize

\subsection{Hyperpfaffian of a sum}

In \cite{Ste}, J. Stembridge gives a formula allowing to compute
the Pfaffian of a sum of  $2s\times 2s$ skew symmetric matrices
\begin{equation}
\Pf(A+B)=\sum_{t=0}^s\sum_{I=(i_1,\cdots,i_{2t})\uparrow}(-1)^{|I|-t}\Pf(A_I
)\Pf(B_{I^c})
\end{equation}
where $|I|=i_1+\cdots i_{2t}$ and $I^c=\{1\cdots 2m\}-I$.\\ This
identity appears to be a special case of a more general one for
hyperpfaffians 
 \begin{equation}\label{stegen}
\Pf^{[2m]}(A+B)=\sum_{k=0}^n\sum_{I=\{i_1<\cdots<i_{2km}\}}(
-1)^{|I|-km} \Pf^{[2m]}(A_I)\Pf^{[2m]}(B_{I^c})
 \end{equation}
 where $A=(a_{i_1,\cdots,i_{2n}})$ and $B=(b_{i_1,\cdots,i_{2n}})$
 are two skew symmetric tensors of order $2m$ over a space of
 dimension $2mn$.\\
 We set
 $$\omega= \sum_{1\leq i_1<\cdots<i_{2m}\leq 2mn}(a_{i_1\cdots
i_{2m}}+b_{i_1\cdots i_{2m}})\eta_{i_1}\cdots\eta_{i_{2m}}\,.$$
The $n{\mbox{th}}$-power of $\omega$ is
$$\omega^n=n!\Pfa^{[2m]}(A+B)\eta_{1}\cdots\eta_{2mn}\,.$$ On the
other hand, we have $$\omega^n=\sum_{k=0}^n\left(n\atop
k\right)\omega_A^k\omega_B^{n-k}\,,$$ where $\omega_A=\sum
a_{i_1\cdots i_{2m}}\eta_{i_1}\cdots\eta_{i_{2m}}$ and
$\omega_B=\sum b_{i_1\cdots
i_{2m}}\eta_{i_1}\cdots\eta_{i_{2m}}$.\\ Since
$$\omega_A^k=k!\sum_{i_1<\cdots<i_{2km}}\Pf^{[2m]}(a_{i_1\cdots
i_{2m}})\eta_{i_1}\cdots\eta_{i_{2km}}\,,$$ we find
$$\omega^n=n!\sum_k\sum_{I=\{i_1<\cdots<i_{2km}\}}\Pf^{[2m]}(A_I)\Pf^{[2m]}(
B_
{I^c})\eta_I \eta_{I^c}\,, $$ where $I^c=\{1,\cdots,2mn\}-I$.
Finally, $$\omega^n=n!\sum_k\sum_{I\uparrow\atop I\subset\{1\cdots
2mn\}}(-1)^{|I|-km} \Pf^{[2m]}(A_I)\Pf^{[2m]}(B_{I^c})\eta_I$$
which by identification leads to (\ref{stegen}).
\subsection{Minor summation formula}
In \cite{IW1}, M. Ishikawa and M. Wakayama  prove a minor
summation formula for Pfaffians
\begin{equation} \sum_{1\leq k_1<\cdots<k_m\leq n}\Pf(B_{k_1\cdots
k_m}^{k_1\cdots k_m})\det\,(T_{k_1\cdots k_m}^{1\cdots
m})=\Pf(Q)\,,\end{equation} where $A^{i_1 \cdots i_m}_{j_1\cdots
j_m }$ denotes the submatrix of $A$ which consists of the rows
$i_1 \cdots i_m$ and the columns $j_1,\cdots, j_m$, $B$ is a skew
symmetric matrix, $T$ is an arbitrary matrix and $Q=TB^tT$. We
shall now extend this formula to hyperpfaffians.

Let $A=(a_{i_1\cdots i_{2m}})_{1\leq i_1,\cdots, i_{2m}\leq 2mn}$
be a skew  symmetric tensor of order $2m$ over a space of
dimension $2mn$ and $T=(T_{ij})_{1\leq i,j\leq 2mn}$ be a square
matrix. We want to find a skew symmetric tensor
$Q=(Q_{i_1,\cdots,i_{2m}})_{1\leq i_1\cdots i_{2n}\leq 2mn}$
verifying
\begin{equation}\label{IWgen}
\sum_{1\leq k_1<\cdots< k_{2mt}\leq 2mn}\Pfa^{[2m]}(A_{k_1\cdots
k_{2mt}})\det(T^{1\cdots 2mt}_{k_1\cdots k_{2mt}})
=\Pf^{[2m]}\left(Q\right)\,,
\end{equation}
where
$A_{k_1,\cdots,k_{2mt}}=(a_{i_1,\cdots,i_{2m}})_{i_1,
\cdots,i_{2m}\in\{k_1,\cdots,k_{2mt}\}}$.\\
In order to compute such a tensor, we need first to remark that
\begin{equation}\label{gendetIW}
\det\left(T_{k_1\cdots k_{2mn}}^{1\cdots
2mn}\right)=\sum_{\sigma\in{\goth
E}_{2mn,2m}}\epsilon(\sigma)\det\left(T^{\sigma(1)\cdots
\sigma(2m)}_{k_1\cdots
k_{2m}}\right)\cdots\det\left(T^{\sigma(2(n-1)m+1)\cdots
\sigma(2mn)}_{k_{2(n-1)m+1}\cdots k_{2mn}}\right)
\end{equation}
which can be obtained by considering
$\omega_k=\sum_{i=1}^{2nt}T_{ik}\eta_i$ and computing
$\omega_{k_1}\cdots\omega_{k_t}$.\\ Now, let us consider
$$Q_{i_1\cdots i_{2m}} =\sum_{1\leq k_1<\cdots<k_{2m}\leq
2mn}a_{k_1\cdots k_{2m}}\det\left(T^{i_1\cdots i_{2m}}_{k_1\cdots
k_{2m}}\right)\,,$$ and the polynomial $$\omega= \sum_{1\leq
i_1<\cdots<i_{2m}\leq 2mt}
\left(\sum_{k_1<\cdots<k_{2m}}a_{k_1\cdots
k_{2m}}\det\left(T^{i_1\cdots i_{2m}}_{k_1\cdots
k_{2m}}\right)\right)\eta_{i_1}\cdots\eta_{i_{2m}}\,.$$ A simple
calculation gives $$w^t=t!\Pf^{[2m]}(Q)\eta_1\cdots\eta_{2mn}$$ On
another hand
\begin{multline}
\omega^t=\sum_{\sigma\in{\goth E}_{2m,2mt}}Q_{i_1\cdots
i_{2m}}\cdots Q_{i_{2m(t-1)+1}\cdots
i_{2mt}}\eta_{i_1}\cdots\eta_{i_{2nt}}\\ =\sum_{1\leq
k_1<\cdots<k_{2mt}\leq 2mn}
 \sum_{\sigma_1\in{\goth E}_{2m,2mt}}a_{k_{\sigma(1)}
\cdots k_{\sigma(2m)}}\cdots a_{k_{\sigma(2m(t-1)+1)}\cdots
a_{\sigma(2mt)}}\\ \sum_{\sigma_2\in{\goth
E}_{2m,2mt}}\epsilon(\sigma_2)\det\left(T^{\sigma_2(1)\cdots
\sigma_2(2m)}_{k_{\sigma_1(1)}\cdots
k_{\sigma_1(2m)}}\right)\cdots\det\left(T^{\sigma_2(2m(t-1)+1)\cdots
\sigma_2(2mt)}_{k_{\sigma_1(2m(t-1)+1)}\cdots
k_{\sigma_1(2mt)}}\right)\\ \times \eta_1\cdots \eta_{2mt}\,.
\end{multline}
By formula (\ref{gendetIW}), we find
\begin{multline}
\omega^t=\sum_{1\leq k_1<\cdots<k_{2mt}\leq
2mn\atop\sigma\in{\goth E}_{2m,2mt}}a_{k_{\sigma(1)} \cdots
k_{\sigma(2m)}}\cdots a_{k_{\sigma(2m(t-1)+1)}\cdots
a_{\sigma(2mt)}} \\ \times \det(T^{1\cdots
2mt}_{k_{\sigma(1)}\cdots
k_{\sigma(2mt)}})\eta_{1}\cdots\eta_{2mt}\,
\end{multline}
so that $$ \omega^t=t!\sum_{1\leq k_1<\cdots< k_{2mt}\leq
2mn}\Pfa^{[2m]}(A_{k_1,\cdots,k_{2mt}})\det(T^{1\cdots
2mt}_{k_1\cdots k_{2mt}})\eta_1\cdots \eta_{2mt}\,,$$ which leads
to the identity (\ref{IWgen}).

\subsection{Generalized de Bruijn formulas}

We consider now a commutative ring $K$ and for $2k$ a non zero
even integer,
 a set of scalars $\{Q_{i_1\cdots i_{2k}}\}_{i_1<\cdots<i_{2k}}$.
We shall first calculate the coefficients $P_I$ of
the $\bigwedge_K(\eta)$-series
\begin{equation}
F=\sum_{I\uparrow}P_I\eta_I=\prod_{i_1}\left(1+\sum_{i_1<
\cdots<
i_k}Q_{i_1 \cdots i_{2k}}\eta_{i_1}\cdots\eta_{i_{2k}}\right).
\end{equation}
As  $\eta_i^2=0$, one has
\begin{equation}
F=\displaystyle\sum_n \sum_{i_{2k(p-1)+1}<\cdots <i_{2kp},\,1\leq
p\leq n\atop i_{2k(q-1)+1}<i_{2kq+1},\,1\leq q\leq n-1}
Q_{i_1\dots i_{2k}}\cdots Q_{i_{2k(n-1)+1}\cdots i_{2kn}}
\eta_{i_1}\cdots \eta_{i_{2nk}}
\end{equation}
and if  $I=\{i_1,\cdots, i_{2nk}\}$ is an increasing sequence, the
coefficient $P_I$ of $\eta_I$ is 
 \begin{equation}\label{xipfa}
P_I={\Pfa}^{[2k]}\left(Q_{i_{j_1}\cdots i_{j_{2k}}}\right) \,.
\end{equation}
Similarly,  in the commutative algebra
$K[\xi]$ ( $\xi^2_i=0$)

\begin{equation}\label{xihaff}
\sum_{I\uparrow\atop l(I)=2nk}
{\Haff}^{[2k]}\left(Q_{i_{j_1}\cdots i_{j_{2k}}}\right)
\xi_I =\prod_{i_1}\left(1+\sum_{i_1<\cdots<i_{2k}} Q_{i_1\cdots
i_{2k}}\xi_{i_1}\cdots\xi_{i_{2k}}\right).
\end{equation}
\bigskip

In the remainder of this section, the ring $K$ is the shuffle
algebra $R\langle A\rangle_{\shuffle}$ for the alphabet
$A=\{a_{i_1 \cdots i_{2k}}\}_{i_1,\cdots,i_{2k}\ge 1}$.We consider
here the series $$
\begin{array}{rcl}
T(A|\eta)&=&\displaystyle\prod_{i_1}^{\rightarrow}\left(1+
\sum_{i_1<\cdots <i_{k~}}\left( \sum_{\sigma\in{\goth
S}_{2k}}\epsilon(\sigma)a_{i_{\sigma(1)}\cdots
i_{\sigma(2k)}}\right) \eta_{i_1}\cdots \eta_{i_{2k}}\right)\\
&=&\displaystyle \sum_{I\uparrow\atop l(I)=2kn}P_I\eta_I
\end{array}$$
 \def\GS{{\bf S}}
 If we set $\GS(i_1,\cdots,i_{2k})=\sum_{\sigma\in{\goth
 S}_{2k}}\epsilon(\sigma)a_{i_{\sigma(1)}\cdots i_{\sigma(2k)}}$,
 we obtain for each increasing sequence $I=(i_1,\cdots,i_{2kn})$
 $$
\begin{array}{rcl}
P_I&=&\displaystyle\int d\eta_{\bar I}T(A|\eta)\\ &=&\displaystyle
\int d\eta_{\bar I}\prod^{\rightarrow}_{i_1\in
I}(1+\sum_{i_2,\cdots,i_{2k}\in I\atop
i_1<\cdots<i_{2k}}\GS(i_1,\cdots,i_{2k})\eta_{i_1}\cdots\eta_{i_{2k}})\\
&=&\displaystyle \sum_{\sigma\in{\goth
E}_{2n,2nk}}\epsilon(\sigma)\GS(i_{\sigma(1)},\cdots,i_{\sigma(2k)})\shuffle
\cdots\shuffle
\GS(i_{\sigma(2(n-1)k+1)},\cdots,i_{\sigma(2nk)})\\ &=&
\displaystyle\sum_{\sigma\in{\goth
S}_{2nk}}\epsilon(\sigma)a_{i_{\sigma(1)}\cdots
i_{\sigma(2k)}}\cdots a_{i_{\sigma(2(n-1)k+1)}\cdots
i_{\sigma(2nk)}}\,.
\end{array}
  $$
  On the other hand, we find by (\ref{xipfa})
  \begin{equation}\label{xipfashu}
P_I={\Pfa}^{[2k]}_{\shuffle} \left(\sum_{\sigma\in{\goth S}_{2k}}
\epsilon(\sigma)a_{{j_{\sigma(1)}} \cdots
i_{{\sigma(2k)}}}\right)_{j_1,\cdots,j_n\in I}
\end{equation}
By Chen's lemma, this equality leads a generalized de Bruijn
formula.
 We consider a set of $4k^2n$
functions $\Phi=\{\phi_{ij}\}_{1\leq i\leq 2k\atop 1\leq j\leq
2kn}$. Equation (\ref{xipfashu}) gives
\begin{equation}\label{Idb1gen}
\begin{array}{rcl}
{\goth W}_{\Phi}(a,b) &=& \displaystyle\mathop
{\int\cdots\int}_{a\leq x_1<\dots<x_n\leq b}
\det(\phi_{1i}(x_j)|\cdots|\phi_{2ki}(x_j))_{1\leq i\leq 2nk\atop
1\leq j\leq n} dx_1\dots dx_n\\
&=&{\Pfa}^{[2k]}\left(\displaystyle \int_a^b\sum_{\sigma\in{\goth
S}_{2k}}\epsilon(\sigma)\phi_{1i_{\sigma(1)}}(x)\cdots
\phi_{2ki_{\sigma(2k)}}(x)dx\right)_{1\leq i_1,\cdots,i_{2k}\leq
2nk }\, .
\end{array}
\end{equation}
Similarly there is a permanent version of this identity
\begin{equation}
\begin{array}{rcl}
{\goth M}_{\Phi}(a,b) &=& \displaystyle\mathop
{\int\cdots\int}_{a\leq x_1<\dots<x_n\leq b} {\rm
per}(\phi_{1i}(x_j)|\cdots|\phi_{2ki}(x_j))_{1\leq i\leq 2nk\atop
1\leq j\leq n} dx_1\dots dx_n\\ &=&{\Hf}^{[2k]}\left(
\displaystyle\int_a^b\sum_{\sigma\in{\goth
S}_{2k}}\phi_{1i_{\sigma(1)}}(x)\cdots
\phi_{2ki_{\sigma(2k)}}(x)dx\right)_{1\leq i_1,\cdots,i_{2k}\leq
2nk }\, .
\end{array}
\end{equation}

\subsection{Sums of even powers of the Vandermonde}
We can use formula (\ref{Idb1gen}) to obtain the average of
certain determinants. For example, let us consider $n$ real
numbers $x_1,\cdots, x_n$ chosen at random from $N$ values
$y_1,\cdots, y_N$ with uniform probability. The average of the
polynomial $\Delta^{2m}(x_1,\cdots,x_n)$ (for $m\in N$) can be
expressed as an integral
\[{\cal
A}(n,m)=\frac1{N^n}\mathop{\int\cdots\int}\Delta^{2m}(x_1
,\cdots,x_n)
\prod_{i=1}^n\left(\sum_{p=1}^N\delta(x_i-y_p)\right)dx_i\] where
$\delta$ denotes the Dirac distribution. We can transform ${\cal
A}(n,m)$ into the interated integral of single determinant
\[{\cal
A}(n,m)=\frac{n!}{N^n}\mathop{\int\cdots\int}_{-\infty<x_1<\cdots<x_n<\infty
}\det\left( f_i^1(x_j)|f_i^2(x_j)|\cdots|f_i^{2m}(x_j)
\right)dx_1\cdots dx_n\] where
\[f_i^1(x)=\left\{\begin{array}{ll}x^{i-1}\sum_{p=1}^N\delta(x-y_p)&\mbox{
if }i\leq n\\ 0&\mbox{ otherwise} \end{array}\right.\] and
\[f_i^s(x)=\left\{\begin{array}{ll}x^{i-(s-1)n-1}&\mbox{ if
}(s-1)n+1\leq i\leq sn\\ 0&\mbox{ otherwise}
\end{array}\right.\] if $1<s\leq 2m$. Using formula
(\ref{Idb1gen}) we find
\begin{equation}\label{average1}
{\cal A}(n,m)=\frac{n!}{N^n}\Pfa^{[2m]}\left(M_{i_1\cdots
i_{2m}}\right)
\end{equation}
where if $i_1<\cdots<i_{2m}$
\begin{equation}
M_{i_1,\cdots,
i_{2m}}=\left\{\begin{array}{ll}\int_{-\infty}^\infty
x^{s(i_1,\cdots,i_{2m})}\sum_{p=1}^N\delta(x-y_p)dx&
\begin{array}{l}\mbox{ if for each } s\leq 2m,\\ (s-1)n<i_s\leq
sn\end{array}\\ 0 &\mbox{ otherwise}\end{array}\right.
\end{equation}
 with
$$s(i_1,\cdots, i_{2m})=\sum_{k=1}^{2m}i_k-m(2n(m-1)+2)$$ and
$$M_{i_{\sigma(1)}\cdots
i_{\sigma(2m)}}=\epsilon(\sigma)M_{i_1\cdots
 i_{2m}}\,.
$$ The nonzero entries of the hyperpfaffian are therefore the
power-sums
\begin{equation}
M_{i_1,\cdots,i_{2m}}=\sum_{p=1}^Ny_p^{s(i_1,\cdots,i_{2m})}.
\end{equation}
with $i_1<\cdots<i_{2m}$ and $(s-1)n<i_s\leq sn$ for
each $s\le 2m$.

In the case where $m=1$, we find
\begin{equation}
{\cal
A}(n,1)=\frac{n!}{N^n}\det\left(\sum_{p=1}^Ny_p^{i+j-2}\right).
\end{equation}
More generally, if we consider the  measure
\[d\mu=\left(\frac
k{N\sqrt\pi}\right)^n\prod_{p=1}^n
\left(\sum_{p=1}^Ne^{-k^2(x_i-y_p)^2}\right)dx_i\,,
\]
we can generalize formula (\ref{average1}) to evaluate
\begin{multline}\label{average1bis}
{\cal A}(k,n,m) =\int\cdots\int\Delta^{2m}(x_1,\cdots, x_n)d\mu\\
=\frac{n!}{N^n}\Pfa^{[2m]}\left(M_{i_1,\cdots,i_{2m}}^{(k)}\right)
\end{multline}
where
\begin{equation}
M_{i_1 \dots
i_{2m}}^{(k)}=\displaystyle\sum_{p=1}^N\sum_{t=0}^{s(i_1,\cdots
,i_{2m})}\left(s(i_1,\cdots,i_{2m})\atop t\right)
y_p^{s(i_1,\dots,i_{2m})-t}\frac{H_t(0)}{(2k)^t}
\end{equation}
where $H_t$ denotes a Hermite polynomial. We can verify that
$\displaystyle{\lim_{k\rightarrow\infty}}{\cal A}(k,n,m)=
{\cal A}(n,m)$, as
expected.

Finally, one can generalize formula (\ref{average1}) by
considering a matrix of functions
\[F=(f_i^k)_{1\leq i\leq 2nm\atop 1\leq k\leq 2m}\]
and the determinant
\[P_F(x_1,\cdots,x_{n})=\det\left(f_i^1(x_j)|\cdots |f_i^{2m}(x_j)\right).\]
Its average, when the  $x_i$ run over the set $\{y_1,\cdots,
y_N\}$ endowed with the uniform probability, can be computed as
above:
\begin{equation}
\begin{array}{rcl}
{\cal
A}(n,m)&=&\displaystyle\frac{n!}{N^n}\mathop{\int\cdots\int}_{-\infty<x_1<
\cdots<x_{2m}<\infty}
P_F(x_1,\cdots,x_{n})\prod_{i=1}^n\left(\sum_{p=1}^N\delta(x_i-y_p)\right)dx
_i\\ &=&\displaystyle\frac{n!}{N^n}
\mathop{\int\cdots\int}_{-\infty<x_1< \cdots<x_{2m}<\infty}
P_{F'}(x_1,\cdots,x_{n})dx_1\cdots dx_n\end{array}
\end{equation}
where \begin{equation} F'(x)=\left(\begin{array}{ccccc}
\sum_{p=1}^N\delta(x-y_p)f^1_1&f_1^2(x)&\cdots&f^{2m}_1(x)\\
\vdots&\vdots&\cdots&\vdots\\
\sum_{p=1}^N\delta(x-y_p)f^1_{2mn}(x)&f^2_{2mn}(x)&\cdots&f^{2m}_{2mn}(x)
\end{array}\right)
\end{equation}
By formula (\ref{Idb1gen}), we find
\begin{equation}
A(n,m) = \frac{n!}{N^n}\Pfa^{[2m]}\left(M_{i_1 \cdots
i_{2m}}\right)
\end{equation}
where
\[
\begin{array}{rcl}
M_{i_1 \cdots i_{2m}}&=&\displaystyle
\int_{-\infty}^\infty\sum_{\sigma\in{\goth S}_{2m}}
\epsilon(\sigma)\sum_{p=1}^N\delta(x-y_p)f^{i_{\sigma(1)}}_1(x)\cdots
f^{i_{\sigma(2m)}}_{2m}(x)dx\\
&=&\displaystyle\sum_{\sigma\in{\goth S}_{2m}}
\epsilon(\sigma)\sum_{p=1}^Nf^{i_{\sigma(1)}}_1(y_p)\cdots
f^{i_{\sigma(2m)}}_{2m}(y_p)
\end{array}
\]

To conclude, let us remark that the hyperpfaffians (\ref{average1})
and (\ref{average1bis}) can alternatively be expressed as hyperdeterminants
in the sense of \cite{Ba}. We note also that in dimension $nk$, the
hyperdeterminant of an alternating tensor $A$ of (even) rank $k$ 
satisfies ${\rm DET}(A)=(\Pfa^{[k]}(A))^k$, generalizing the 
classical case $k=2$. Such relations will be investigated in
a forthcoming paper.

\bigskip\small
\sc
\noindent Institut Gaspard Monge\\
Universit\'e de Marne-la-Vall\'ee\\
5 Boulevard Descartes\\
Champs-sur-Marne\\
77454 Marne-la-Vall\'ee cedex 2\\
FRANCE\\

\end{document}